\documentclass[a4paper]{article}
\usepackage[turkish]{babel}
\usepackage[latin5]{inputenc}
\usepackage[T1]{fontenc}
\usepackage{amsthm}
\usepackage{amsfonts}
\usepackage{mathrsfs}

\theoremstyle{definition}
\newtheorem{teo}{Theorem}
\newtheorem{lemma}{Lemma}
\newtheorem{ex}{Example}
\begin{document}

\title{The Bounds for Eigenvalues of Normalized and
Signless Laplacian Matrices }
\author{ŞERİFE BÜYÜKKÖSE AND ŞEHRİ GÜLÇİÇEK ESKİ}
\date{June,2014}
\maketitle

In this paper, we obtain the bounds of the extreme
eigenvalues of a normalized and signless Laplacian
matrices using by their traces. In addition, we determine the bounds for k-th
eigenvalues of normalized and signless Laplacian matrices.
\vspace{11pt} \newline

\section{Introduction}

Let $G(V,E)$ be a simple graph with the vertex set $V=\{v_1,v_2,\ldots,v_n\}$
and edge set of $E$. For $v_i \in V$, the degree of $v_i$, the set of
neighbours of $v_i$ denoted by $d_i$ and $N_i$, respectively. The
cardinality of $N_i$ is denoted by $c_{ij}$; i.e, $|N_i|=c_{ij}$. If $v_i$
and $v_j$ adjacency, we denote $v_i\sim v_j$ of shortly use $i\sim j$.

The adjacency matrix, Laplacian matrix and diagonal matrix of vertex degree
of a $G$ graph denoted by $A(G)$, $L(G)$, $D(G)$, respectively. Clearly
\[
L(G)=D(G)-A(G).
\]

The normalized Laplacian matrix of $G$ is defined as $\mathscr{L}%
(G)=D^{-1/2}(G)L(G)D^{-1/2}(G)$ i.e, $\mathscr{L}(G)=[\ell_{ij}]_{n\times n}$,
where
\[
\ell_{ij}=\left\{
\begin{array}{cll}
1 & ; & if \ \ i=j \\
\frac{-1}{\sqrt{d_{i}d_{j}}} & ; & if \ \ i\sim j \\
\ 0 & ; & otherwise.%
\end{array}
\right.
\]
\newline

The signless Laplacian matrix of $G$ is defined as $Q(G)=D(G)+A(G)$ i.e, $%
Q(G)=[q_{ij}]_{n\times n}$, where
\[
q_{ij}=\left\{
\begin{array}{cll}
d_{i} & ; & if \ \ i=j \\
1 & ; & if \ \ i\sim j \\
\ 0 & ; & otherwise.%
\end{array}
\right.
\]

Since $\mathscr{L}(G)$ normalized Laplacian matrix and $Q(G)$ signless
Laplacian matrix are real symetric matrices, their eigenvalues are real. We
denote the eigenvalues of $\mathscr{L}(G)$ and $Q(G)$ are by
\[
\lambda_{1}(\mathscr{L}(G)) \geq \cdots \geq \lambda_{n}(\mathscr{L}(G))
\]
and
\[
\lambda_{1}(Q(G)) \geq \cdots \geq \lambda_{n}(Q(G))
\]
,respectively.

Now we give some bounds for normalized Laplacian matrix and signless
Laplacian matrix.\newline

\textbf{Oliveira and de Lima's bound [1].} For a simple connected
graph $G$ with $n$ vertices and $m$ edges, $\Delta=d_1\geq d_2\geq\cdots
\geq d_n=\delta$
\begin{eqnarray}
\lambda_1(Q(G))\leq \max_i \left\{\frac{d_i+\sqrt{d_i^2+8d_i m_i}}{2}\right\}
\end{eqnarray}
where $m_i=\frac{1}{d_i}\sum_{i\sim j}d_j$.

\textbf{Another Oliveira and de Lima's bound [1].}
\begin{eqnarray}
\lambda_1(Q(G))\leq \max_i \left\{d_i+\sqrt{d_i m_i} \right\}
\end{eqnarray}
where $m_i=\frac{1}{d_i}\sum_{i\sim j}d_j$.

\textbf{Li, Liu et al. bound's [2,3].}

\begin{eqnarray}
\lambda_1(Q(G))\leq \frac{\Delta+\delta-1\sqrt{(\Delta+%
\delta-1)^2+8(2m-(n-1)\delta)}}{2}.
\end{eqnarray}

\textbf{Rojo and Soto's bound [4].} If $\lambda_1$ is the largest
eigenvalue of $\mathscr{L}$, then
\begin{eqnarray}
|\lambda_{1}(\mathscr{L}(G))|\leq 2-\min_{i<j} \left(\frac{|N_i\cap N_j|}{%
\max\{d_i,d_j\}} \right)
\end{eqnarray}
where the minimum is taken over all pairs $(i,j)$, $(1\leq i<j\leq n)$.

In this paper, we find an extreme eigenvalues of normalized Laplacian
matrix and signless Laplacian matrix of a $G$ graph with using their traces.

To obtain bounds for eigenvalues of $\mathscr{L}(G)$ and $Q(G)$ we need the
following lemmas and theorems.

\begin{lemma}
Let $W$ and $\lambda=(\lambda_{j})$ be nonzero column vectors, $%
e=(1,1,\ldots ,1)^T$, $C=I_n-\frac{ee^T}{n}$ $m=\frac{\lambda^{T}e}{n}$, $%
s^2=\frac{\lambda^TC\lambda}{n}$ and $I_n$ is an Identity matrix. Let $%
\lambda_{1}\geq\lambda_{2}\geq\cdots\geq\lambda_{n}$. Then
\[
-s\sqrt{nW^{T}CW}\leq W^{T}\lambda-mW^{T}e=W^{T}C\lambda\leq s\sqrt{nW^{T}CW}%
.
\]
\[
\sum_{j}(\lambda_j-\lambda_n)^{2}=n[s^{2}+(m-\lambda_n)^{2}]
\]
\[
\sum_{j}(\lambda_1-\lambda_j)^{2}=n[s^{2}+(\lambda_1- m)^{2}].
\]
\[
\lambda_{n}\leq m-\frac{s}{\sqrt{n-1}}\leq m+\frac{s}{\sqrt{n-1}}\leq
\lambda_{1}.
\]
\end{lemma}

\begin{teo}
Let A be a $n\times n$ complex matrix. Conjugate transpose of $A$ denoted by
$A^{\ast}$. Let $B=AA^{\ast}$ whose eigenvalues are $\lambda_{1}(B)\geq%
\lambda_{2}(B)\geq\cdots\geq\lambda_{n}(B).$ Then
\[
m-s\sqrt{n-1}\leq \lambda_{n}^{2}(B)\leq m-\frac{s}{\sqrt{n-1}}
\]
and
\[
m+\frac{s}{\sqrt{n-1}}\leq\lambda_{1}^{2}(B)\leq m+s\sqrt{n-1}
\]
where $m=\frac{trB}{n}$ \ and \ $s^{2}=\frac{trB^{2}}{n}-m.$
\end{teo}

\section{Main Results for Normalized Laplacian Matrix}

\begin{teo}
Let G be a simple graph and $\mathscr{L}(G)$ be a normalized Laplacian
matrix of G. If the eigenvalues of $\mathscr{L}(G)$ are $\lambda_{1}(%
\mathscr{L}(G)) \geq \lambda_{2}(\mathscr{L}(G)) \geq \cdots \geq
\lambda_{n}(\mathscr{L}(G))$, then

\begin{eqnarray}
\lambda_{n}(\mathscr{L}(G)) \leq \sqrt{\left(1+\frac{2}{n}\sum_{i\sim j ,
i<j}\frac{1}{d_{i}d_{j}}\right)+\sqrt{\frac{tr[L(G)]^{4}-nm^{2}}{n(n-1)}}}
\end{eqnarray}
\begin{eqnarray}
\lambda_{1}(\mathscr{L}(G))\geq \sqrt{\left(1+\frac{2}{n}\sum_{i\sim j , i<j}%
\frac{1}{d_{i}d_{j}}\right)+\sqrt{\frac{tr[L(G)]^{4}-nm^{2}}{n(n-1)}}}
\end{eqnarray}

\begin{eqnarray}
\lambda_{1}(\mathscr{L}(G)) \leq \sqrt{1+\frac{2}{n}\sum_{i\sim j , i<j}%
\frac{1}{d_{i}d_{j}}+\sqrt{\left(\frac{tr[L(G)]^{4}}{n}-m^{2}\right)(n-1)}}
\end{eqnarray}
\end{teo}

\emph{Proof.} Obviously,
\[
tr[\mathscr{L}(G)]^{2}=n+2\sum_{i\sim j , i<j}\frac{1}{d_{i}d_{j}}
\]
and
\[
tr[\mathscr{L}(G)]^{4}=\sum_{i=1}^{n} \left(1+ \sum_{i\sim j}\frac{1}{d_i d_j%
} \right)^2 + 2\sum_{i<j}\left( \sum_{k \in N_{i}\cap N_{j}} \frac{1}{d_k%
\sqrt{d_id_j}}- \sum_{i\sim j}\frac{2}{\sqrt{d_id_j}} \right)^2
\]
Since $\mathscr{L}(G)$ real symmetric matrix, we find the result from
Theorem 1.

\begin{ex}
Let $G=(V,E)$ with $V=\{1,2,3,4,5,6\}$ and \newline
$E=\{(1,2),(1,5),(2,3),(2,4),(2,6),(3,4),(3,5),(4,5),(5,6)\}.$

\begin{tabular}{l|c|c|c}
$\lambda(\mathscr{L}(G))$ & (4) & (6)(lower bound) & (7)(upper bound) \\
\hline
1.86 & 2 & 1.34 & 1.93 \\
&  &  &
\end{tabular}
\end{ex}

\section{Main Results for Signless Laplacian Matrix}

\begin{teo}
Let G be a simple graph and $Q(G)$ be a signless Laplacian matrix of G. If
the eigenvalues of $Q(G)$ are $\lambda_{1}(Q(G)) \geq \lambda_{2}(Q(G)) \geq
\cdots \geq \lambda_{n}(Q(G))$, then
\begin{eqnarray}
\lambda_{n}(Q(G)) \leq \sqrt{\left(1+\frac{2}{n}\sum_{i\sim j , i<j}\frac{1}{%
d_{i}d_{j}}\right)+\sqrt{\frac{tr[Q(G)]^{4}-nm^{2}}{n(n-1)}}}
\end{eqnarray}
\begin{eqnarray}
\lambda_{1}(Q(G))\geq \sqrt{\left(1+\frac{2}{n}\sum_{i\sim j , i<j}\frac{1}{%
d_{i}d_{j}}\right)+\sqrt{\frac{tr[Q(G)]^{4}-nm^{2}}{n(n-1)}}}
\end{eqnarray}

\begin{eqnarray}
\lambda_{1}(Q(G)) \leq \sqrt{1+\frac{2}{n}\sum_{i\sim j , i<j}\frac{1}{%
d_{i}d_{j}}+\sqrt{\left(\frac{tr[Q(G)]^{4}}{n}-m^{2}\right)(n-1)}}
\end{eqnarray}
\end{teo}

\emph{Proof.} Clearly
\[
tr[Q(G)]^{2}=n+2\sum_{i\sim j , i<j}\frac{1}{d_{i}d_{j}}
\]
and
\[
tr[Q(G)]^{4}=\sum_{i=1}^{n} \left(1+ \sum_{i\sim j}\frac{1}{d_i d_j}
\right)^2 + 2\sum_{i<j}\left( \sum_{k \in N_{i}\cap N_{j}} \frac{1}{d_k\sqrt{%
d_id_j}}- \sum_{i\sim j}\frac{2}{\sqrt{d_id_j}} \right)^2
\]
Since $Q(G)$ real symmetric matrix, we found the result from Theorem 1.

\begin{ex}
Let $G=(V,E)$ with $V=\{1,2,3,4,5,6,7\}$ and \newline
$E=\{(1,2),(1,3),(1,4),(1,5),(1,6),(1,7),(2,3),(3,5),(4,5),(4,6)\}.$

\begin{tabular}{l|c|c|c|c|c}
$\lambda(Q(G))$ & (1) & (2) & (3) & (9)(lower bound) & (10)(upper bound) \\
\hline
7.67 & 9.08 & 9.74 & 9.34 & 4.58 & 7.76 \\
&  &  &  &  &
\end{tabular}
\end{ex}

\ \

\begin{center}
\textbf{REFERENCES}
\end{center}

\begin{enumerate}
\item[{\textbf{[1]}}] C.S. Oliveira, L.S.De Lima, N.M.M. De Abreu, P. Hansen,
\emph{Bound on the Index of the Signless Laplacian of a graph, Discrete
Appl. Math.} \textbf{158}, (2010), 355-360

\item[{\textbf{[2]}}] J. Li, Y. Pan, \emph{Upper Bounds for the Laplacian
Graph \newline
Eigenvalues, Acta Math. Sin. Engl. ser. 20(5)}, (2004), 803-806

\item[{\textbf{[3]}}] H. Liu, M. Lu, F. Tian, \emph{On the Laplacian Spectral
Radius of a Graph, Linear Algebra Appl.}, (2004) ,\textbf{376}, 135-141

\item[{\textbf{[4]}}] O. Rojo, R.L. Soto, \emph{A New Upper Bound On The
Largest \newline
Normalized Laplacian Eigenvals, Operators and Matrices}, (2013) number2,
volume7, 323-332

\item[{\textbf{[5]}}] H. Wolkowich, G.P.H. Styan, \emph{Bounds for
Eigenvalues Using Traces of Matrice}, \emph{Linear Algebra and its
Applications}, (1980), \textbf{29}, 471-506
\end{enumerate}

\begin{flushright}
Şerife BÜYÜKKÖSE \newline
Gazi University Faculty of Science \newline
Department of Mathematics \newline
06500 Ankara-Turkey\newline
e-mail: sbuyukkose@gazi.edu.tr\newline
serifebuyukkose@gmail.com\newline
\end{flushright}

\begin{flushright}         
Şehri Gülçiçek ESKİ\newline
Ahi Evran University\newline
Institue of Science \newline
40100 Kırşehir-Turkey\newline
e-mail: gulcicekeski@gmail.com\newline
\end{flushright}

\end{document}